\documentclass[11pt]{article}

\usepackage[english]{style}
\usepackage{latexsym,amsmath,amssymb,amsfonts,diagrams,graphics,amscd}
\usepackage{graphicx}
\input xypic
\xyoption{all}
\input epsf.tex

\usepackage{pgf,tikz}
\usetikzlibrary{arrows}
\usepackage{amssymb}
\DeclareFontFamily{U}{rsf}{}
\DeclareFontShape{U}{rsf}{m}{n}{
  <5> <6> rsfs5 <7> <8> <9> rsfs7 <10-> rsfs10}{}
\DeclareMathAlphabet{\mathscr}{U}{rsf}{m}{n}

\DeclareMathAlphabet{\mathgth}{U}{euf}{m}{n}

\DeclareFontFamily{U}{cyr}{}
\DeclareFontShape{U}{cyr}{m}{n}{
  <5> wncyr5 <6> wncyr6 <7> wncyr7 <8> wncyr8 <9> wncyr9 <10-> wncyr10}{}
\DeclareMathAlphabet{\mathcyr}{U}{cyr}{m}{n}

\input cyracc.def

\makeatletter
\def\operator@font{\sf}
\makeatother

\newcommand{\sbt}{{\scalebox{0.5}{\textbullet}}}
\newcommand{\sA}{{\mathcal A}}
\newcommand{\sB}{{\mathcal B}}

\newcommand{\sD}{{\mathcal D}}

\newcommand{\sL}{{\mathcal L}}

\newcommand{\cI}{{\mathscr I}}

\newcommand{\cO}{{\mathscr O}}

\newcommand{\bW}{{\overline{W}}}
\newcommand{\bS}{{\overline{S}}}

\newcommand{\bi}{{\bar{i}}}

\newcommand{\ccH}{{\mathscr H}}

\newcommand{\dR}{{\mathsf dR}}

\newcommand{\sHom}{\underline{\mathsf{Hom}}}
\newcommand{\RsHom}{{\R\sHom}}

\newcommand{\sEnd}{\underline{\mathsf{End}}}

\newcommand{\D}{{\mathbf D}}

\newcommand{\R}{\mathbf{R}}

\renewcommand{\S}{{\mathbb{S}}}

\DeclareMathOperator{\Spec}{Spec}

\newcommand{\Hom}{{\mathsf{Hom}}}

\DeclareMathOperator{\Tor}{Tor}

\DeclareMathOperator{\Pic}{Pic}

\DeclareMathOperator{\id}{id}

\DeclareMathOperator{\Ext}{Ext}

\newcommand{\ra}{\rightarrow}

\newcommand{\inj}{\hookrightarrow}

\newcommand{\C}{\mathbb{C}}

\newcommand{\Z}{\mathbf{Z}}

\newcommand{\iso}{\cong}

\newcommand{\del}{\partial}

\newcommand{\field}[1]{\mathbb{#1}}
\renewcommand{\P}{\field{P}}
\newcommand{\A}{\field{A}}

\renewcommand{\phi}{\varphi}

\newarrow{Equal} =====

\title{Hodge theorem for the logarithmic de Rham complex via derived intersections}
\author{M\'arton Hablicsek\thanks{Mathematics Department, University of Pennsylvania,
David Rittenhouse Lab, 209 S.\ 33rd Street, Philadelphia, PA 19104,
USA, {\em e-mail: }{\tt mhabli@math.upenn.edu}}}
\begin{document}
\maketitle
In a beautiful paper \cite{DelIll} Deligne and Illusie proved the degeneration of the Hodge-to-de Rham spectral sequence using positive characteristic methods.  In a recent paper \cite{AriCalHab2} Arinkin, C\u ald\u araru and the author of this paper gave a geometric interpretation of the problem of Deligne-Illusie showing that the triviality of a certain line bundle on a derived scheme implies the the Deligne-Illusie result.  In the present paper we generalize these ideas to logarithmic schemes and using the theory of twisted derived intersection of logarithmic schemes we obtain the Hodge theorem for the logarithmic de Rham complex.
\section{Introduction}

\paragraph
Let $Y$ be a smooth proper variety over an algebraically closed field $k$ of characteristic 0.  The algebraic de Rham complex is defined as the complex
\[\Omega^\sbt_Y:=0\ra \cO_Y\xrightarrow{d} \Omega^1_{Y/k}\xrightarrow{d}...\]
where $d$ is the usual differential on the {\em algebraic} forms.  The de Rham cohomology of $Y$ is defined as the hypercohomology of the de Rham complex,
\[H^*_{\dR}(Y)=R^*\Gamma(Y,\Omega^\sbt_Y).\]
The Hodge-to-de Rham spectral sequence
\[{}^1E_{pq}=H^p(Y,\Omega^q_Y)\Rightarrow H^{p+q}_{\dR}(Y)\]
is given by the stupid filtration on the de Rham complex whose associated graded terms are the $\Omega^q_Y$.

\paragraph
In their celebrated paper \cite{DelIll}, Deligne and Illusie proved the degeneration of the Hodge-to-de Rham spectral sequence holds in positive characteristics by showing the following result.

\begin{Theorem}\cite{DelIll}
Let $X$ be a smooth proper scheme over a perfect field $k$ of
positive characteristic $p > \dim X$. Assume that $X$ lifts to the ring $W_2(k)$ of
second Witt vectors of $k$. Then the Hodge-to-de Rham spectral sequence
for $X$ degenerates at ${}^1E$.
\end{Theorem}
\medskip

Then they showed that the corresponding result in characteristic 0 follows from a standard reduction argument.

\paragraph
Similar results can be obtained in the logarithmic setting.  Consider a smooth proper scheme $X$ over a perfect field $k$ of characteristic $p>\dim X$ and a reduced normal crossing divisor $D$ on $X$.  We can find local coordinates $x_1,...,x_n$ around each point of the divisor so that in an \'etale neighborhood of that point the divisor is cut out by $x_1\cdot...\cdot x_k$ for some $k\leq n$.  The logarithmic 1-forms are generated locally by the symbols
\[\frac{dx_1}{x_1},...,\frac{dx_k}{x_k},dx_{k+1},...,dx_n\]
and the sheaf of logarithmic 1-forms is denoted by $\Omega^1_X(\log D)$.  The sheaf $\Omega^q_X(\log D)$ of logarithmic $q$-forms is defined as $\wedge^q\Omega_X^1(\log D)$.  The differential $d$ of the {\em meromorphic} de Rham complex maps logarithmic forms to logarithmic forms, and hence we define the logarithmic de Rham complex $\Omega^\sbt_X(\log D)$ as the subcomplex of the meromorphic de Rham complex consisting of logarithmic forms.  The stupid filtration of this complex gives rise to a Hodge-to-de Rham spectral sequence
\[{}^1E_{pq}=H^p(Y,\Omega^q_X(\log D))\Rightarrow R^{p+q}\Gamma(X,\Omega^\sbt_X(\log D)).\]

In \cite{Kat}, Kato generalized the result of Deligne-Illusie obtaining the following result.

\begin{Theorem}\cite{Kat}
\label{thm:Kat}
Assume that the pair $(X,D)$ lifts to the ring $W_2(k)$.  Then the Hodge-to-de Rham spectral sequence of the logarithmic de Rham complex degenerates at ${}^1E$.
\end{Theorem}
\medskip

As before the corresponding result in characteristic 0 follows from a standard reduction argument.  (For a detailed treatment, and for applications to vanishing theorems see \cite{EsnVie}.)

\paragraph In a recent paper \cite{AriCalHab2}, Arinkin, C\u ald\u araru and the author of the present paper recast the problem of the degeneration of the Hodge-to-de Rham spectral sequence as a derived self-intersection problem.  They considered the embedding of the Frobenius twist $X'$ of $X$ into its cotangent bundle $T^*X'$ as the zero section.  The sheaf of crystalline differential operators $D_X$ on $X$ can be regarded as a sheaf $D$ of Azumaya algebras on $T^*X'$.  This Azumaya algebra splits on the zero section $X'\ra T^*X'$ giving rise to an embedding of Azumaya spaces $X'\ra (T^*X',D)$.  The intersection problem given by a smooth subvariety $Y$ inside an Azumaya space $\bS=(S,\sA)$ is called a twisted intersection problem in \cite{AriCalHab2}.  The main geometric observation in \cite{AriCalHab2} is that there exists a line bundle on the ordinary derived self-intersection of $Y$ inside $S$ measuring the difference between the twisted and the ordinary derived self-intersections.  In the case of $X'\ra (T^*X',D)$ this line bundle is given by the dual of $F_*\Omega^\sbt_X$, the Frobenius pushforward of the de Rham complex of $X$.  

Using the theory of twisted derived self-intersections the authors obtain the following result giving a geometric interpretation of the result of Deligne-Illusie (see also \cite{OguVol} for the relation between the statements (1) and (2)).

\begin{Theorem}\cite{AriCalHab2}
\label{thm:AriCalHab2} Let $X$ be a smooth scheme over a perfect field $k$ of characteristic $p > \dim X$. Then the following five statements are equivalent.
\begin{itemize}
\item[(1)] $X$ lifts to $W_2(k)$.
\item[(2)] $D$ splits on the first infinitesimal neighborhood of $X'$ in $T^*X'$.
\item[(3)] The associated line bundle is trivial.
\item[(4)] $F_*\Omega^\sbt_X$ is formal in $D(X')$ (meaning that it is quasi-isomorphic to the direct sum of its cohomology sheaves).
\end{itemize}
\end{Theorem}
\medskip

\paragraph In the present paper we generalize the result of \cite{AriCalHab2} to the logarithmic setting.  We consider the logarithmic scheme $(X,D)$ where $X$ is a smooth proper scheme over a perfect field $k$ of characteristic $p>\dim X$ and $D$ is a reduced normal crossing divisor on $X$.  As the category of quasi-coherent sheaves on $(X,D)$ we consider the quasi-coherent parabolic sheaves on $(X,D)$ (see \cite{Yok}).  Parabolic sheaves were first introduced by Mehta and Seshadri (\cite{MehSes}, \cite{Ses}) on a projective curve with finitely many marked points as locally-free sheaves $E$ with filtrations
\[0=F_k(E)\inj F_{k-1}(E)\inj...\inj F_0(E)=E\]
at every marked point in order to generalize to non-projective curves the correspondence between stable sheaves of rank 0 and irreducible unitary representations of the topological fundamental group.

\paragraph
The sheaf of crystalline logarithmic differential operators $D_X(\log D)$ is not an Azumaya algebra over its center.  On the other hand one can equip $D_X(\log D)$ with a filtration so that the corresponding parabolic sheaf is a parabolic Azumaya algebra over its center.  (For a treatment of parabolic Azumaya algebras see \cite{KulLie}.)  The center of the Azumaya algebra can be identified with the structure sheaf of the logarithmic cotangent bundle $T^*X'(\log D')$ of the Frobenius twist $(X',D')$ of $(X,D)$ equipped with the divisor $\pi^*D'$ where $\pi:T^*X'(\log D')\ra X'$ is the bundle map.  Therefore the parabolic sheaf of crystalline logarithmic differential operators can be regarded as a parabolic sheaf $D(\log D)_*$ over $(T^*X'(\log D'),\pi^*D')$.  As before we consider the embeddings 
\begin{itemize}
\item $(X',D')\ra (T^*X'(\log D'),\pi^*D')$ and 
\item $(X',D')\ra (T^*X'(\log D'),\pi^*D',D(\log D)_*)$.
\end{itemize}
The difference between the derived self-intersections of the embeddings is measured by a parabolic line bundle, which is the dual of the Frobenius pushforward of the logarithmic de Rham complex $F_*\Omega^\sbt_X(\log D)$ equipped with the filtration given by $D$.

We generalize the theory of twisted derived self-intersections to the logarithmic setting.  Our main result is the generalization of Theorem \ref{thm:AriCalHab2}.

\begin{Theorem}
\label{thm:equi1}
Let $X$ be a smooth variety over a perfect field of characteristic $p>\dim X$, with a reduced normal crossing divisor $D$.  Then, the following statements are equivalent.
\begin{itemize}
\item[(1)] The logarithmic scheme $(X,D)$ lifts to $W_2(k)$.
\item[(2)] The associated line bundle is trivial.
\item[(3)] The parabolic sheaf of algebras $D(\log D)_*$ splits on the first infinitesimal neighborhood of $(X',D')$ inside $(T^*X'(\log D'),\pi^*D')$.
\item[(4)] The complex $F_*\Omega^\sbt_X(\log D)_*$ is quasi-isomorphic to a formal parabolic sheaf equipped with the trivial parabolic structure.
\end{itemize}
\end{Theorem}
\medskip

\noindent
As an easy corollary we obtain Theorem \ref{thm:Kat}.

\paragraph
\label{par:Konque} The main application we have mind is answering a question of Kontsevich.  Let $X\ra \P^1$ be a rational function, so that $f^{-1}(\A^1)$ is a smooth complex algebraic variety and $f^{-1}(\infty)$ is a normal crossing divisor.  In the above setting Kontsevich introduced a family of complexes $(\Omega^\sbt_X,ud+v\wedge df)$, $(u,v)\in \C^2$ generalizing the notion of a twisted de Rham complex.  He conjectured that the hypercohomology spaces are independent of the choices of $u$ and $v$.  In a recent preprint \cite{KatKonPan} Katzarkov, Kontsevich and Pantev verify this conjecture in the case of a reduced normal crossing divisor $D$, and in \cite{EsnSabYu} a complete proof is given for any normal crossing divisor.  We believe that our methods can be generalized to give another proof of this question of Kontsevich.

\paragraph
We remark that there is another approach in the literature to deal with the sheaf of differential operators in the logarithmic setting.  In \cite{Sch} Schepler extends the theory of Ogus and Vologodsky (\cite{OguVol}) to the case of logarithmic schemes.  He uses the theory of indexed modules and algebras and he shows the the sheaf of differential operators form an indexed Azumaya algebra.  These results generalized by Ohkawa (\cite{Ohk1}, \cite{Ohk2}) to the ring of differential operators for higher level.  We choose to not work with indexed modules and algebras because of lack of functoriality.  It would be interesting to understand the results of \cite{Sch} , \cite{Ohk1} and \cite{Ohk2} from our approach.

\paragraph The paper is organized as follows.  In Section \ref{sec:log} we collect some basic facts about logarithmic schemes and parabolic sheaves.  In Section \ref{sec:charp} we introduce the parabolic sheaf of crystalline differential operators and the parabolic logarithmic de Rham complex.  We show the Azumaya property of the parabolic sheaf of crystalline differential operators and that there exists a Koszul duality between the parabolic sheaf of crystalline differential operators and the parabolic logarithmic de Rham complex.  In Section \ref{sec:derint} we summarize the theory of twisted derived intersections, which we briefly expand to the case of logarithmic schemes in Section \ref{sec:logderint}.  We conclude the paper with Section \ref{sec:mainthm} where we prove our main theorems Theorem \ref{thm:equi1} and Theorem \ref{thm:Kat}.

\paragraph 
\textbf{Acknowledgements.} The author expresses his thanks to Dima Arinkin, Andrei C\u ald\u araru and Tony Pantev for useful conversations.
\section{Background on logarithmic schemes}
\label{sec:log}

In this section we investigate the notion of the category of (quasi-)coherent sheaves on logarithmic schemes $(X,D)$: the category of (quasi-)coherent parabolic sheaves.  We remark that we could have taken an alternative path, looking at (quasi-)coherent sheaves on the infinite root stack (see \cite{BorVis}, \cite{TalVis}).  We follow the notations of \cite{Yok}.

\paragraph We regard $\R$ as the category whose objects are real numbers and whose morphism spaces between two objects $\alpha\in \R$ and $\beta\in \R$ are defined as
\[Mor(\alpha,\beta)=\begin{cases} \{i^{\alpha,\beta}\} & \mbox{if }\alpha\geq \beta,\\ \emptyset & \mbox{otherwise}\end{cases}.\]

\begin{Definition}
An $\R$-filtered $\cO_X$-module is a covariant functor from the category $\R$ to the category of $\cO_X$-modules.  For an $\R$-filtered $\cO_X$-module $E$ we denote the $\cO_X$-module $E(\alpha)$ by $E_\alpha$, and the $\cO_X$-linear homomorphisms $E(i^{\alpha,\beta})$ by $i_E^{\alpha,\beta}$.  For an $\R$-filtered $\cO_X$-module $E$ we define the $\R$-filtered $\cO_X$-module $E[\alpha]$ as $E[\alpha]_{\beta}=E_{\alpha+\beta}$ with morphisms $i^{\beta,\gamma}_{E[\alpha]}=i_E^{\beta+\alpha,\gamma+\alpha}$.  In the sequel we denote $\R$-filtered $\cO_X$-modules by $E_*$.
\end{Definition}

\begin{Definition}
For an $\R$-filtered $\cO_X$-module $E_*$ and for an "ordinary" $\cO_X$-module $F$ we define their tensor product as $(E_*\otimes F)_\alpha=E_\alpha\otimes F$ with homomorphisms $i^{\alpha,\beta}_{E_*\otimes F}=i^{\alpha,\beta}_{E_*}\otimes \id_F$.
\end{Definition}
\medskip

\noindent
We are ready to define the category of (quasi-)coherent parabolic sheaves with respect to an effective Cartier divisor $D$.

\begin{Definition}
A (quasi-)coherent {\em parabolic} sheaf is an $\R$-filtered $\cO_X$-module $E_*$ together with an isomorphism of $\R$-filtered $\cO_X$-modules
\[E_*\otimes \cO_X(-D)\iso E_*[1].\]
Parabolic morphisms between parabolic sheaves $E_*$ and $E'_*$ are natural transformations $E_*\ra E'_*$.  We denote the set of parabolic morphisms by $Hom(E_*, E'_*)$.
\end{Definition}
\medskip 

\paragraph[Remark:]\label{rmk:01} It is enought to know $E_\alpha$ and $i_E^{\alpha,\beta}$ for $\alpha,\beta\in [0,1]$ to determine the parabolic sheaf $E_*$.  We say that a parabolic sheaf has weights $\alpha=(\alpha_0,...\alpha_k)$
\[0=\alpha_0<\alpha_1<...<\alpha_k<1\]
if $E_\beta=E\gamma$ and $i_E^{\beta,\gamma}=\id$ for all $\beta, \gamma\in [0,1]$ satisfying $\alpha_i< \beta,\gamma\leq \alpha_{i+1}$.

For us, there are natural choices for $k$ and the $\alpha_i$, and thus, we define the category of (quasi-)coherent sheaves on $(X,D)$ as follows.

\begin{Definition} The category of (quasi-)coherent sheaves on $(X,D)$ is the category whose objects are parabolic sheaves with $k=p$ and weights $\alpha_i=\frac{i}{p+1}$, and the Hom-spaces are the sets of parabolic morphisms.  These categories will be denoted by $Coh(X,D)$ or $QCoh(X,D)$ respectively.
\end{Definition}

\paragraph The category $QCoh(X,D)$ is abelian, the kernel and cokernel of a morphism can be defined pointwise.  It has enough injectives \cite{Yok}, we denote the corresponding derived category by $\D(X,D)$.  

\paragraph
Consider a morphism $f:X\ra Y$ and an effective Cartier divisor $D$ on $Y$ such that $f^*D$ is an effective divisor on $X$.  This data gives rise to a morphism of logarithmic schemes $ (X,f^*D)\ra (Y,D)$.  We abuse notation and denote the induced map by $f$ as well.  We define the pushforward and pullback along $f$ as follows.  For any parabolic sheaf $E_*$ on $X$ its pushforward $f_*E_*$ is defined as the parabolic sheaf where $(f_*E)_\alpha$ are the pushforward of the sheaves $E_\alpha$ along $f$ and the morphisms $i^{\alpha,\beta}_{f_*E}$ are the morphisms $f_*i^{\alpha,\beta}_E$.  Indeed we obtain a parabolic sheaf on $(Y,D)$, by the projection formula, we have
\[f_*(E_*[1])=f_*(E_*\otimes \cO_X(-f^*D))=f_*E_*\otimes \cO_Y(-D)=f_*E_*[1].\]
Similarly, for any parabolic sheaf $E'_*$ on $Y$ its pullback $f^*E'_*$, is defined as the parabolic sheaf where $(f^*E')_\alpha$ is the pullback of $E'_\alpha$ along $f$ and the morphisms $i^{\alpha, \beta}_{f^*E'}$ are the morphisms $f^*i^{\alpha,\beta}_{E'}$.  Again, we obtain a parabolic sheaf on $(X,f^*D)$, we have
\[f^*(E'_*[1])=f^*(E'_*\otimes \cO_Y(-D))=f^*E'_*\otimes \cO_X(-f^*D)=f^*E'_*[1].\]
The pushforward and pullback functors descend to the derived categories and by abuse of notation we denote the corresponding maps by $f_*$ and $f^*$ as well.

\paragraph[Remark:]In general, the pushforward and the pullback morphisms of a parabolic sheaf along a morphism $f:(X,D_1)\ra (Y,D_2)$ are more complicated than as above.  Our case is special, we have $f^*D_2=D_1$. 

\paragraph The categories $Coh(X,D)$ and $QCoh(X,D)$ are equipped with natural monoidal structures.  In order to define the monoidal structures we take a quick detour.

An important subcategory of $QCoh(X,D)$ is the category of parabolic bundles (\cite{Ses}, \cite{MehSes}).

\begin{Definition}
A {\em parabolic bundle} is a triple $(E,F_*,\alpha_*)$ where $E$ is a locally-free sheaf on $X$, $F_*$ is a filtration of $E$ by coherent sheaves on $X$
\[F_k(E)=E\otimes \cO_X(-D)\inj F_{k-1}(E)\inj F_{k-2}(E)\inj...\inj F_0(E)=E\]
together with a sequence of weights $\alpha$ satisfying
\[0=\alpha_0<\alpha_1<...<\alpha_k<1.\]
\end{Definition}
\medskip

\noindent
The sequence of weights determines a family of coherent sheaves $E_x$ for $0\leq x\leq 1$ defined as
\[E_0=E\quad\mbox{and}\quad E_x=F_i(E)\]
for $\alpha_i< x\leq \alpha_{i+1}$.  A morphism between parabolic bundles $(E,F_*,\alpha_*)$ and $(E',F'_*,\alpha'_*)$ is a morphism of $\cO_X$-modules $\phi:E\ra E'$ so that $\phi(E_x)\subseteq E'_x$ for any $x\in [0,1]$.  By Remark \ref{rmk:01} parabolic bundles give rise to parabolic sheaves and morphisms between parabolic bundles are exactly the parabolic morphisms between the corresponding parabolic sheaves.

\paragraph
Consider the morphism 
\[\psi:\Pic X\times \Z\left[\frac{1}{p}\right]\ra Coh(X,D)\]
mapping the pair $(L,a)$ to the parabolic bundle $(L,F_*,\alpha_*)$ where
\[L_x=\begin{cases}L & \mbox{if }x\leq a',\\ L\otimes \cO_X(-D) & \mbox{if }a'<x\leq 1.\end{cases}\]
Here $a'$ denotes the residue of $a$ modulo 1.  Parabolic bundles of this form are the parabolic line bundles.  The tensor product of parabolic sheaves is defined (for parabolic line bundles) to respect the group structure coming from the natural group structure on $\Pic(X)\times \Z\left[\frac{1}{p}\right]$ given by
\[(L_1,a_1)\cdot (L_2,a_2)=(L_1\otimes L_2,a_1+a_2).\]
The unit element of the tensor product is given by the parabolic sheaf $\psi(\cO_X,0)$, where
\[\psi(\cO_X,0)_x=\begin{cases}\cO_X & \mbox{if }x=0,\\ \cO_X(-D) & \mbox{if }0<x\leq 1.\end{cases}\]
We define the {\em structure sheaf} of the logarithmic scheme $(X,D)$ to be the parabolic sheaf $\psi(\cO_X,0)$ and in the sequel we denote it by $\cO_{(X,D)}$.

\begin{Definition}
\label{log:hom} For two parabolic sheaves $E_*$, $F_*$ we define the sheaf Hom functor as 
\[\sHom_x(E_*,F_*):=\Hom(E_*,F_*[x]).\]
\end{Definition}

In particular, for any parabolic line bundle $L$, its parabolic sheaf of endomorphisms is isomorphic to $\cO_{(X,D)}$.  The sheaf Hom functor and the tensor product satisfy the usual adjoint property (for more details, see \cite{Yok}) giving rise to natural monoidal structures on $Coh(X,D)$ and $QCoh(X,D)$.  We remark that the monoidal structure descends to the derived category $\D(X,D)$.
\section{Background on schemes over fields of positive characteristics}
\label{sec:charp}

In this section we collect basic facts about schemes over fields of positive characteristics.  We review the notion of logarithmic tangent sheaf, logarithmic $q$-forms and the crystalline sheaf of logarithmic differential operators.  We show that the parabolic sheaf of crystalline logarithmic differential operators is an Azumaya algebra over its center.  We conclude the section by showing that there exists a Koszul duality between the parabolic sheaf of crystalline logarithmic differential operators and the parabolic logarithmic de Rham complex.

\paragraph Let $X$ be a smooth scheme over a perfect field $k$ of characteristic $p$.  The absolute Frobenius map $\phi:\Spec k\ra \Spec k$ is associated to the $p$-th power map $k\ra k$.  The Frobenius twist of $X$ is defined as the base change of $X$ along the absolute Frobenius morphism
\[\xymatrix{X'\ar[r]\ar[d]& X\ar[d]\\ \Spec k\ar[r]^\phi&\Spec k.}\]
The $p$-th power map $\cO_X\ra \cO_X$ gives rise to a morphism $X\ra X$ compatible with $\phi$ and thus it factors through the Frobenius twist.  The induced morphism $F:X\ra X'$ is called the relative Frobenius morphism.  For any effective Cartier divisor, $D$ on $X$ with ideal sheaf $\cI=\cO_X(-D)$; we obtain a corresponding Cartier divisor $D'$, which is the pullback of $D$ under the base change morphism $X'\ra X$ and whose pullback $F^*D'$ under the relative Frobenius morphism is the divisor $pD$.  Thus, we have the following sequence of maps of logarithmic schemes
\[(X,D)\ra (X,pD)\ra (X',D')\ra (X,D).\]
The morphism $(X,D)\ra (X',D')$ is called the relative Frobenius morphism of logarithmic schemes, by abuse of notation, we denote it by $F$ as well.

\paragraph  We say that a derivation $\delta\in T_X$ is {\em logarithmic} if for every open subset $U$ we have $\delta(\cI(U))\subset \cI(U)$.  The logarithmic derivations form a subsheaf of the tangent bundle of $X$ which is called the logarithmic tangent sheaf, $T_X(\log D)$.  The sheaf $T_X(\log D)$ is a Lie subalgebroid of $T_X$ meaning that it is closed under the Lie-bracket on $T_X$.  In characteristic $p>0$ the $p$-th iteration $\delta^{[p]}$ of a derivation $\delta$ is again a derivation.  Clearly, the sheaf $T_X(\log D)$ is closed under this operation as well making $T_X(\log D)$ a sub-$p$-restricted Lie algebroid of $T_X$.

In general $T_X(\log D)$ is not a subbundle of $T_X$, we say that a divisor is {\em free} if $T_X(\log D)$ is a locally free sheaf.  For instance, reduced normal crossing divisors are free divisors.  Indeed, we can find local coordinates $x_1,...,x_n$ around each point of the divisor so that in an \'etale neighborhood of that point the divisor is given by $x_1\cdot ...\cdot x_k$ for some $k\leq n$, and thus the logarithmic tangent sheaf is generated by the logarithmic derivations $x_1\frac{\del}{\del x_1}, x_2\frac{\del}{\del x_2},...,x_k\frac{\del}{\del x_k}$ and $\frac{\del}{\del x_{k+1}},...,\frac{\del}{\del x_n}$.

\paragraph[Remark:]  Consider the parabolic bundle $\sL_{(X,pD)}$ on the logarithmic scheme $(X,pD)$ defined as the parabolic sheaf $(\cO_X,F_*,\alpha_*)$ where the filtration is given by the natural filtration 
\[\cO_X(-pD)\inj \cO_X(-(p-1)D)\inj...\inj \cO_X(-D)\inj \cO_X.\]
We remark that this parabolic bundle is the pushforward of $\cO_{(X,D)}$ under the natural morphism $(X,D)\ra (X,pD)$.  The corresponding parabolic bundle on $(X',D')$ is the parabolic bundle $(F_*\cO_X,F_*,\alpha_*)$ where the filtration is given by
\[F_*\cO_X\otimes \cO_{X'}(-D')=F_*\cO_X(-pD)\inj...\inj F_*\cO_X(-D)\inj F_*\cO_X.\]
This parabolic bundle is the pushforward of $\cO_{(X,D)}$ along the relative Frobenius morphism $F:(X,D)\ra (X',D')$, and hence we denote it by $F_*\cO_{(X,D)}$.

It is easy to see that those derivations of $\cO_X$ which respect the filtration of the parabolic bundle $\sL_{(X,pD)}$ are exactly the logarithmic derivations of $X$.

\paragraph We say that a meromorphic $q$-form $\omega$ is logarithmic, if for every affine open subset where $\cI=(g)$ for some $g\in \cO_X$, we have that $\omega\wedge dg$ and $d(\omega)$ are algebraic.  The logarithmic $q$-forms form a subsheaf $\Omega^q_X(\log D)$ of the sheaf of {\em meromorphic} $q$-forms.  The differential $d$ of the meromorphic de Rham complex maps logarithmic forms to logarithmic forms.  We define the logarithmic de Rham complex $\Omega^\sbt_X(\log D)$ as the subcomplex of the de Rham complex consisting of the logarithmic forms.  This complex $\Omega^\sbt_X(\log D)$ is not a complex of $\cO_X$-modules, the differential $d$ is not linear in $\cO_X$.  On the other hand, we have
\[d(s^p\omega)=ps^{p-1}ds\wedge \omega+s^pd\omega=s^pd\omega\]
for every $s\in \cO_X$ and $\omega\in \Omega^q_X(\log D)$.  This implies that $F_*\Omega^\sbt_X(\log D)$ is a complex of $\cO_{X'}$-modules.  If the divisor is free, then the sheaves $\Omega^q_X(\log D)$ are locally free sheaves and moreover we have $\Omega^q_X(\log D)=\wedge^q \Omega^1_X(\log D)$.  In the case of a reduced normal crossing divisor, locally the logarithmic $1$-forms are generated by $d(\log x_1),...,d(\log x_k)$ and $dx_{k+1},...,dx_{n}$.  Similarly to the non-logarithmic case, there is a perfect duality between $T_X(\log D)$ and $\Omega^1_X(\log D)$ given by contracting with polyvector fields.  

In the sequel $D$ denotes a reduced normal crossing divisor.

\paragraph
\label{par:center} The sheaf of crystalline logarithmic differential operators $D_X(\log D)$ is defined as the universal enveloping algebra of the Lie algebroid $T_X(\log D)$.  Explicitly it is defined locally as the $k$-algebra generated by sections of $T_X(\log D)$ and $\cO_X$ modulo the relations
\begin{itemize}
\item $s\cdot\delta=s\delta$ for every $s\in \cO_X$ and $\delta\in T_X(\log D)$,
\item $\delta_1\cdot \delta_2-\delta_2\cdot\delta_1=[\delta_1,\delta_2]$ for every $\delta_1, \delta_2\in T_X(\log D)$ and
\item $\delta\cdot s-s\cdot\delta=\delta(s)$ for every $s\in \cO_X$ and $\delta\in T_X(\log D)$. 
\end{itemize}
We emphasize that we do not work with the sheaf of PD differential operators, for our purposes we need an algebra which is of finite type over $X$.  Since $T_X(\log D)$ is a Lie subalgebroid of $T_X$, we have an inclusion of $D_X(\log D)$ into the sheaf of crystalline differential operators $D_X$ (defined as the universal enveloping algebra of $T_X$).  

\paragraph The map
\[\psi: T_X(\log D)\ra D_X(\log D)\]
mapping
\[\delta\mapsto \delta^p-\delta^{[p]}\]
is $\cO_{X'}$-linear and its image is in the center of $D_X(\log D)$ (see \cite{BezMirRum} for a detailed treatment in the non-logarithmic case) implying that the center of $D_X(\log D)$ can be identified with the structure sheaf  $\cO_{T^*X'(\log D')}$ of the logarithmic cotangent bundle $\pi: T^*X'(\log D')\ra X'$ over the Frobenius twist.  The zero section $i:X'\ra T^*X'(\log D')$ of the bundle map $\pi$ gives rise to a natural embedding of logarithmic schemes
\[i_D:(X',D')\ra (T^*X'(\log D'),\pi^*D'),\]
since $i^*\pi^*D'=D'$.

\paragraph We equip the sheaf of algebras $D_X(\log D)$ with the trivial logarithmic structure on $(X,D)$, we denote the corresponding logarithmic sheaf by $D_X(\log D)_*$:
\[D_X(\log D)_x=\begin{cases} D_X(\log D) & \mbox{if }x=0,\\D_X(\log D)\otimes \cO_X(-D) &\mbox{if }0<x\leq 1.\end{cases}\]

After pushing forward the parabolic sheaf $D_X(\log D)_*$ along the relative Frobenius morphism $F:(X,D)\ra (X',D')$ we obtain the parabolic sheaf $F_*D_X(\log D)_*$, whose filtration is given by
\[F_*D_X(\log D)\otimes \cO_{X'}(-D')=F_*(D_X(\log D)\otimes \cO_X(-pD))\inj...\inj F_*(D_X(\log D)).\]
This parabolic sheaf has weights $\alpha_i=\frac{i}{p+1}$, similarly to the weights of $\sL_{(X,pD)}$ or of $F_*\cO_{(X,D)}$.

\paragraph
By the discussion in Paragraph \ref{par:center} we can regard the parabolic sheaf of algebras $F_*D_X(\log D)_*$ as a parabolic sheaf of algebras on $(T^*X'(\log D'), \pi^*D')$.  We denote the corresponding parabolic sheaf of algebras by $D(\log D)_*$.  The bundle map $\pi$ identifies $\pi_*D(\log D)_*$ with $F_*D_X(\log D)_*$.

The following lemma is a straightforward generalization of the non-logarithmic result in \cite{BezMirRum}.

\begin{Lemma}
\label{lem:az}
Assume that $D$ is a reduced normal crossing divisor.  Then, the parabolic sheaf of algebras $D(\log D)_*$ is an Azumaya algebra over the logarithmic space $(T^*X'(\log D'),\pi^*D')$.  Moreover, the Azumaya algebra is split on the zero section $(X',D')\ra (T^*X'(\log D'),\pi^*D')$.
\end{Lemma}
\medskip

\begin{Proof}
We only highlight the key steps.  We need to show that $D(\log D)_*$ becomes a matrix algebra under a flat cover.  Consider an affine open set $U'$ of $X'$, and the corresponding open set $U\subseteq X$.  Pick local coordinates $x_1,...,x_n$ for $U$ so that the reduced normal crossing divisor, $D$ is given by the equation $x_1\cdot...\cdot x_k$.  Then, the (non-parabolic) sheaf of algebras $D_X(\log D)$ is generated by $\Gamma(U,\cO_X)$ and the derivations $x_i\frac{d}{dx_i}$ for $1\leq i\leq k$, and the derivations $\frac{d}{dx_i}$ for $k<i\leq n$.  Consider $F_*D_X(\log D)$ and the centralizer $A_X$ of $\Gamma(U',\cO_{X'})$ inside $F_*D_X(\log D)$.  A straightforward calculation shows that the centralizer is the $R=\Gamma(U',\cO_{X'})$-algebra generated by the logarithmic derivations $\delta$ of $X$, in other words
\[A_X=R\left[x_1\frac{d}{dx_1},...,x_k\frac{d}{dx_k},\frac{d}{dx_{k+1}},...,\frac{d}{dx_n}\right]\]
There is a natural logarithmic structure on $V=\Spec A_X$ given by the restriction of the divisor $\pi^*D'$ to $V$.  The ideal sheaf corresponding to the divisor is generated by the element $x_1^p\cdot...\cdot x_k^p$.  We denote by $\cI$ the ideal sheaf of $\cO_X$ generated by the element $x_1\cdot...\cdot x_k$.  The sheaf 
\[D(\log D)|_V:=D(\log D)\otimes_{\cO_{T^*X'(\log D')(U)}} A_X\]
is generated by $\Gamma(U',\cO_{X'})$ and two copies $u_1,...,u_n$ and $v_1,...,v_n$ corresponding to the logarithmic derivations
\[x_1\frac{d}{dx_1},...,x_k\frac{d}{dx_k},\frac{d}{dx_{k+1}},...,\frac{d}{dx_n}.\] 
The action $u_i.x=u_ix$ and $v_i.x=xv_i$ for $x\in D(\log D)$ gives rise to an action of $D(\log D)|_V$ on the parabolic sheaf
\[\sD:=D(\log D)\otimes \cI^p\inj D(\log D)\otimes \cI^{p-1}\inj...\inj D(\log D)\]
viewed as a parabolic sheaf over $\Spec A_X$.  A local calculation similar to one in \cite{BezMirRum} shows that we have an isomorphism
\[D(\log D)_x|_V=\sEnd_{(V,\pi^*D'|_V)}(\sD)_x\]
hence $D(\log D)_*$ is an Azumaya algebra over $(T^*X'(\log D'),\pi^*D')$.
\medskip 

Next we show that $D(\log D)_*|_{X'}$ is a split Azumaya algebra.  More explicitly, we show that $D(\log D)_*|_{X'}=\sEnd_{(X',D')}(F_*\cO_{(X,D)})$.  We remind the reader that the parabolic sheaf $F_*\cO_{(X,pD)}$ is defined as the parabolic sheaf given by the filtration
\[F_*\cO_X\otimes \cO_{X'}(-D')=F_*\cO_X(-pD)\inj...\inj F_*\cO_X(-D)\inj F_*\cO_X.\]
The sheaf of crystalline logarithmic differential operators acts non-trivially on $\cO_X(-D)$, $\cO_X(-2D)$, ..., $\cO_X(-(p-1)D)$ and moreover the elements which act trivially are generated by the symbols $\delta^p-\delta^{[p]}$ for $\delta\in T_X(\log D)$.  These are exactly the elements which vanish under the pullback map $\cO_{T^*X'(\log D')}\ra \cO_{X'}$, hence $D(\log D)_*|_{X'}$ acts on $F_*\cO_{(X,D)}$.  As before a local calculation similar to one in \cite{BezMirRum} shows that we have an isomorphism
\[D(\log D)_*|_{X'}=\sEnd_{(X',D')}(F_*\cO_{(X,D)}).\]
This concludes the proof.\qed
\end{Proof}

\paragraph[Remark:] If $D$ is a non-reduced normal crossing divisor, then the lemma above does not hold anymore. It would be interesting to generalize the above lemma to any normal crossing divisor (see Paragraph \ref{par:Konque}).

\paragraph As a consequence of Lemma \ref{lem:az}, we obtain a Morita equivalence between the category of coherent sheaves $Coh(X',D')$ on $(X',D')$ and the category of coherent sheaves $Coh(D(\log D)|_{(X',D')})$ on $(X',D')$ with a left $D(\log D)|_{X'}$-action: the functors 
\[m_*:Coh(X',D')\ra Coh(D(\log D)|_{(X',D')}):\]
\[E_*\mapsto E_*\otimes F_*\cO_{(X,D)}\]
and
\[m^*:Coh(D(\log D)|_{(X',D')})\ra Coh(X',D'):\]
\[ E_*\mapsto \sHom_{D(\log D)|_{X'}}(F_*\cO_{(X,D)},E_*)\]
are inverses to each other.  These functors give rise to an equivalence between the corresponding derived categories $\D(X',D')$ and $\D(X', D', D(\log D)|_{(X',D')})$.

\paragraph We conclude this section by showing that the sheaf of crystalline logarithmic differential operators and the logarithmic de Rham complex are Koszul dual (our reference is \cite{CalNar}).  The sheaf $\cO_X$ is naturally a left $D_X(\log D)\subset D_X$-module given by the action of the logarithmic derivations on $\cO_X$.  The logarithmic Spencer complex $Sp^\sbt(\cO_X)$ defined as the complex of left $D_X(\log D)$-modules
\[0\ra D_X(\log D)\otimes \wedge^n T_X(\log D)\ra ...\ra D_X(\log D)\otimes T_X(\log D)\ra D_X(\log D)\]
where the differentials \
\[d_{Sp}:D_X(\log D)\otimes \wedge^i T_X(\log D)\ra D_X(\log D)\otimes \wedge^{i-1} T_X(\log D)\]
are given by
\begin{align*}d_{Sp}(T\otimes \delta_1\wedge ...\wedge \delta_i)&=\sum_{l=1}^i (-1)^{i-1} T\delta_l \otimes \delta_1\wedge...\wedge \hat{\delta_l}\wedge ...\wedge \delta_i+\\
&+\sum_{l<k=1}^i (-1)^{l+k} T\otimes [\delta_l,\delta_k]\wedge \delta_1\wedge...\wedge\hat{\delta_l}\wedge...\wedge\hat{\delta_k}\wedge ...\wedge \delta_i\end{align*}
is a locally free resolution of $\cO_X$ by locally free left $D_X(\log D)$-modules.  As a consequence, for any left $D_X(\log D)$ module $F$ we have that the object $\RsHom_{(X,D_X(\log D))} (\cO_X,F)$ can be represented by the logarithmic de Rham complex of $F$
\[0\ra F\ra F\otimes \Omega^1_X(\log D)\ra...\ra F\otimes \Omega^n_X(\log D).\]
In the case of $F=\cO_X$ we obtain an isomorphism 
\[\Omega^\sbt_X(\log D)=\RsHom_{(X,D_X(\log D))} (\cO_X,\cO_X).\]
Similarly, for any  left $D_X(\log D)$-module $E$ we obtain a complex $Sp^\sbt(E)$ defined as
\[0\ra D_X(\log D)\otimes \wedge^n T_X(\log D)\otimes E\ra ...\ra D_X(\log D)\otimes E\]
where the differentials 
\[d:D_X(\log D)\otimes \wedge^i T_X(\log D)\otimes E\ra D_X(\log D)\otimes \wedge^{i-1} T_X(\log D)\otimes E\]
are given by
\begin{align*}d(T\otimes \delta_1\wedge ...\wedge \delta_i\otimes e)&=d_{Sp}(T\otimes\delta_1\wedge...\wedge \delta_i)\otimes e-\\
&-\sum_{l=1}^i (-1)^{i-1} T \otimes \delta_1\wedge...\wedge \hat{\delta_l}\wedge ...\wedge \delta_i\otimes \delta_l(e).\\\end{align*}
\paragraph The Spencer complex of $E$ gives rise to a locally free resolution of $E$ by left $D_X(\log D)$-modules.  As a consequence we can compute 
\[\RsHom_{(X,D_X(\log D))}(\cO_X(lD),\cO_X(mD))\]
for any $l,m\in \Z$ by replacing $\cO_X(lD)$ by its Spencer complex, and we obtain isomorphisms
\begin{equation}\label{eq:iso}\RsHom_{(X,D_X(\log D))}(\cO_X(lD),\cO_X(mD))=\Omega^\sbt_X(\log D)\otimes \cO_X((m-l)D).\end{equation}

\paragraph
The above discussion shows that for the parabolic sheaf $D_X(\log D)$ we have
\[\RsHom_{(X,D,D_X(\log D)_*)}(\cO_{(X,D)},\cO_{(X,D)})_*=\Omega^\sbt_X(\log D)_*\]
where $\Omega^\sbt_X(\log D)_*$ is the parabolic sheaf given by the de Rham complex equipped with the trivial parabolic structure:
\[\Omega^\sbt_X(\log D)_x=\begin{cases} \Omega^\sbt_X(\log D) & \mbox{if }x=0,\\ \Omega^\sbt_X(\log D)\otimes \cO_X(-D) & \mbox{if }0<x\leq 1.\end{cases}\]
Consider the parabolic sheaf $F_*D_X(\log D)_*$ on $(X',D')$.  We remind the reader that this parabolic sheaf has weights $\alpha_i=\frac{i}{p+1}$.  As before, using the isomorphisms \ref{eq:iso} we obtain an isomorphim
\[\RsHom_{(X',D',F_*D_X(\log D)_*)}(F_*\cO_{(X,D)},F_*\cO_{(X,D)})_*=F_*\Omega^\sbt_X(\log D)_*\]
where $F_*\Omega^\sbt_X(\log D)_*$ is the complex of parabolic bundles given by the pushforward de Rham complex with the filtration
\[F_*\Omega^\sbt_X(\log D)\otimes \cO_{X'}(-D')\inj ...\inj F_*(\Omega^\sbt_X(\log D)\otimes \cO_X(-D))\inj F_*\Omega^\sbt_X(\log D).\]
Similarly to the parabolic sheaf $F_*\cO_{(X,D)}$, the parabolic sheaf $F_*\Omega^\sbt_X(\log D)_*$ has weights $\alpha_i=\frac{i}{p+1}$.

\section{Derived self-intersection of (Azumaya) schemes}
\label{sec:derint}

In this section we summarize the theory of derived self-intersections of (Azumaya) schemes.  Our references are \cite{AriCal}, \cite{AriCalHab1}, \cite{CioKap} and \cite{Gri}.

\paragraph Let $S$ be a smooth variety and $X$ be a smooth subvariety of $S$ of codimension $n$.  Assume that the base field is either of characteristic 0 or of $p>n$.  We denote the embedding of $X$ inside $S$ by $i$ and the corresponding normal bundle by $N$.  The derived self-intersection $W$ of $X$ inside $S$ is a dg-scheme whose structure sheaf is constructed by taking the derived tensor product of the structure sheaf of $X$ with itself over $\cO_S$.  The derived self-intersection is equipped with a map from the underived self-intersection, $X$.  

\paragraph We say that a derived scheme $W$ is {\em formal} over a scheme $\pi: W\ra Z$ if $\pi_*\cO_W$ is a {\em formal} complex of $\cO_Z$-modules, meaning that there exists an isomorphism of commutative differential graded algebras
\[\pi_*\cO_W=\bigoplus_k \ccH^k(\pi_*\cO_W)[-k].\] 
A local calculation \cite{CalKatSha} shows that the cohomology sheaves of the structure sheaf of the derived self-intersection $W$ (over $S$) are given by 
\[\ccH^{-*}(\cO_W)=\Tor_*^S(\cO_X,\cO_X)=\wedge^* N^\vee.\]
Therefore, the formality of the derived self-intersection asserts that there is a quasi-isomorphism (of commutative dg-algebras)
\[\pi_*\cO_W=\bigoplus_{i=0}^n \wedge^i N^\vee[i]=:\S(N^\vee[1]).\]
(We omit writing the pushforward of $N^\vee$ to $Z$ along the map $X\ra W\ra Z$.)
The main result of \cite{AriCal} is the following.

\begin{Theorem}[\cite{AriCal}]
\label{thm:Wformal}
The following statements are equivalent.
\begin{itemize}
\item[(1)] There exists an isomorphism of dg-autofunctors of $D(X)$
\[i^*i_*(-)=(-)\otimes \S(N^\vee[1]).\]
\item[(2)] $W$ is formal over $X\times X$.
\item[(3)] The natural map $X\ra W$ is split over $X\times X$. 
\item[(4)] The short exact sequence 
\[0\ra T_X\ra T_S|_X\ra N\ra 0\]
of vector bundles on $X$ splits.
\end{itemize}
\end{Theorem}
\medskip

\paragraph[Remark:] For instance, the derived self-intersection of $X'$ inside $T^*X'(\log D')$ is formal: the bundle map $\pi:T^*X'(\log D')\ra X'$ splits the embedding of the zero section, and hence the injection 
\[T_{X'}\inj T_{T^*X'(\log D')}|_{X'}\]
is split.

\paragraph An {\em Azumaya scheme} is a pair $(S,\sA)$ where $S$ is a scheme or dg-scheme and $\sA$ is an Azumaya algebra over $S$.  We say that $(f,E):(X,\sB)\ra (S,\sA)$ is a 1-morphism of Azumaya schemes if $f$ is a morphism of schemes and $E$ is an $f^*\sA^{opp}\otimes \sB$-module that provides a Morita equivalence between $f^*\sA$ and $\sB$.  Given an embedding of Azumaya schemes $\bi:(X,\sA|_X)\ra (S,\sA)$, the derived self-intersection is given by $\bW=(W,\sA|_W)$.

\paragraph 
From now on, assume that $\sA|_X$ is a split Azumaya algebra with splitting module $E$, in other words there exists a 1-isomorphism $(\id, E)$ of Azumaya schemes $m:X\ra (X,\sA|_X)$.  We denote the induced map $X\ra (S,\sA)$ by $i'$.  We organize our spaces into the following diagram.
\[\xymatrix{W\ar[rr]\ar[dd]&&X\ar[d]^m\ar@/^2pc/[dd]^{i'}\\&(W,\sA|_W)\ar[r]^p\ar[d]^q& (X,\sA|_X)\ar[d]^\bi\\X\ar[r]^-m\ar@/_2pc/[rr]^{i'}&(X,\sA|_X)\ar[r]^{\bi}& (S,\sA)}\] 
\medskip

\paragraph The Azumaya schemes $W$ and $(W,\sA|_W)$ are abstractly isomorphic, but in general the isomorphism is not over $(X,X)$.  (We remark that $W$ can be thought of as a dg-scheme over $X\times X$, on the other hand Azumaya spaces do not have absolute products, and thus it is more natural to think of them as spaces equipped with two morphisms to $X$.)  The structure sheaves of derived self-intersections $W$ and $(W,\sA|_W)$ regarded as dg-schemes endowed with a map to $(X,X)$ are the kernels of the dg-autofunctors of $D(X)$ $i^*i_*(-)$ and $i'^*i'_*(-)$ respectively.  In \cite{AriCalHab1} the authors show that there exists an isomorphism of dg-autofunctors of $D(X)$
\[i'^*i'_*(-)=q_*p^*(-\otimes L)\]
for some line bundle $L$ on the derived scheme $W$.  This line bundle is called the associated line bundle of the derived self-intersection in \cite{AriCalHab1}. 

\paragraph Consider the object $\bi^*\bi_*E$ of $D(X,\sA|_X)$.  A local calculation similar to one in \cite{CalKatSha} shows that there exist isomorphisms
\[\ccH^k(\bi^*\bi_*E)=E\otimes \wedge^{-k} N^\vee[-k].\]
Therefore, we obtain a triangle in $D(X,\sA|_X)$
\[E\otimes N^\vee[1]\ra \tau^{\geq -1}\bi^*\bi_*E\ra E\ra E\otimes N^\vee[2].\]
The rightmost map of the triangle
\[\alpha_E\in H^2_{(X,\sA|_X)}(E,E\otimes N^\vee)=H^2(X,N^\vee)\]
is called the HKR class of $E$.  

\paragraph
The HKR class gives the obstruction of lifting $E$ to the first infinitesimal neighborhood of $(X,\sA|_X)$ inside $(S,\sA)$.  Having a lifting of $E$ is equivalent to having a splitting module $F$ of the Azumaya algebra $\sA|_{X^{(1)}}$ so that $F|_X=E$ where $X^{(1)}$ denotes the first infinitesimal neighborhood of $X$ inside $S$.  Such lifting exists, if $i'$ splits to first order meaning that there exists a map $\phi:X^{(1)}\ra X$ splitting the natural inclusion $X\ra X^{(1)}$, so that $\phi^*E$ is a splitting module for $\sA|_{X^{(1)}}$.  We are ready to state the main result concerning about the triviality of the associated line bundle.

\begin{Theorem}[\cite{AriCalHab1}]
\label{thm:Azu}
Assume that $W$ is formal over $X\times X$.  Then, the following statements are equivalent.
\begin{itemize}
\item[(1)] The dg-schemes $W$ and $(W,\sA|_W)$ are isomorphic over $X\times X$.
\item[(2)] There exists an isomorphism of dg-autofunctors of $D(X)$
\[i^*i_*(-)\iso i'^*i'_*(-)=(-)\otimes \S(N^\vee[1]).\]
\item[(3)] The associated line bundle is trivial.
\item[(4)] The morphism $i'$ splits to first order.
\item[(5)] The HKR class $\alpha_E$ vanishes.
\end{itemize}
\end{Theorem}
\medskip

\section{Derived self-intersection of logarithmic (Azumaya) schemes}
\label{sec:logderint}

In this section we expand the theory of twisted derived intersections to the logarithmic setting.  

\paragraph Let $X$ be a smooth subscheme of a smooth scheme $S$.  We denote the embedding by $i$.  Let us equip $S$ with an effective Cartier divisor $D$ so that $i^*D$ is an effective divisor on $X$.  Consider the induced embedding of logarithmic schemes $i_D:(X,D|_X)\ra (S,D)$.  We say that $i_D$ splits to first order if there is a left inverse of the induced morphism
\[(X,D|_X)\ra (X^{(1)},D|_{X^{(1)}}),\]
where $X^{(1)}$ denotes the first infinitesimal neighborhood of $X$ inside $S$.  Equivalently, we say that $i_D$ splits to first order if there is a splitting $\rho:X^{(1)}\ra X$ of the embedding $X\ra X^{(1)}$ so that $\rho^*D|_X=D|_{X^{(1)}}$. 

\paragraph
We define the derived-self intersection of $(X,D|_X)$ inside $(S,D)$ as the logarithmic dg-scheme $(W,D|_W)$.  We assume that the dg-scheme $W$ is formal over $X\times X$.  Notice that the divisor $D|_W$ can be thought of restricting $D|_{X\times X}$ to $W$.  We generalize Theorem \ref{thm:Wformal} to the logarithmic setting.

\begin{Proposition}
\label{prop:for}
Assume further that $i_D$ splits to first order.  Then we have an isomorphism of dg-autofunctors of $D(X,i^*D)$
\[i_D^*i_{D,*}(-)=(-)\otimes \S(N^\vee[1]).\]
In other words, $(W,D|_W)$ is formal over $(X\times X,D|_{X\times X})$.
\end{Proposition}
\medskip 

\begin{Proof}
The proof is entirely similar to of the proof of Theorem 0.7 of \cite{AriCal}.\qed
\end{Proof}

\paragraph[Remark:]\label{rem:ass} All the assumptions of the Proposition above are satisfied in the case when $S$ is a vector bundle $\pi:S\ra X$ over $X$, $i:X\ra S$ is the zero section and the divisor $D$ is the pullback along $\pi$ of a divisor on $X$.

\paragraph
\label{par:setup} We turn our attention to embeddings of logarithmic Azumaya schemes.  Let us equip the logarithmic scheme $(S,D)$ with a parabolic sheaf of Azumaya algebras $\sA_*$ and assume that $\sA_*$ splits over $(X,D|_X)$ with splitting module $E_*$.  We remind the reader that $E_*$ induces an isomorphism of spaces $m_D:(X,D)\ra(X,D|_X,\sA_*|_X)$.  We denote the embedding of logarithmic Azumaya spaces $(X,D|_X,\sA_*|_X)\ra (S,D,\sA_*)$ by $\bar{i}_D$, and the composite of embeddings $(X,D|_X)\ra (X,D|_X,\sA_*|_X)\ra (S,D,\sA_*)$ by $i_D'$.  We organize our spaces as follows.
\[\xymatrix{(W,D|_W)\ar[rr]^{p_D}\ar[dd]_{q_D}&&(X,D)\ar[d]^{m_D}\ar@/^3pc/[dd]^{i'_D}\\&(W,D|_W,\sA_*|_W)\ar[r]\ar[d]& (X,D|_X,\sA_*|_X)\ar[d]^{\bi_D}\\(X,D)\ar[r]^-{m_D}\ar@/_3pc/[rr]^{i'_D}&(X,D|_X,\sA_*|_X)\ar[r]^{\bi_D}& (S,D,\sA_*)}\] 
\medskip

\paragraph The spaces $(W,D|_W)$ and $(W,D|_W,\sA_*|_W)$ are abstractly isomorphic, but in general not over the pair $((X,D), (X,D))$, since the splitting modules $p_D^*E_*$ and $q_D^*E_*$ may not be isomorphic.  Two splitting modules differ by a parabolic line bundle, thus the failure to have an isomorphism between $(W,D|_W)$ and $(W,D|_W,\sA_*|_W)$ is measured by a parabolic line bundle on $(W,D|_W)$ which we call the associated parabolic line bundle $\sL_*$.  As a consequence, we obtain an isomorphism of dg-functors
\[i_D'^*i'_{D,*}(-)=q_{D,*}p_D^*(-\otimes \sL_*).\]
In particular for the structure sheaf $\cO_{(X,D|_X)}$ we have
\[i_D'^*i'_{D,*}\cO_{(X,D|_X)}=q_{D,*}p_D^*\sL_*.\]

\paragraph Consider the object $\bar{i}_D^*\bar{i}_{D,*}E_*$.  As before, a local calculation similar to one in \cite{CalKatSha} shows that there exist isomorphisms of parabolic sheaves
\[\ccH^{k}(\bi_D^*\bi_{D,*}E_*)=E_*\otimes \wedge^{-k} N^\vee[-k].\]
As above, we define the HKR class $\alpha^D_{E_*}$ as the rightmost map of the triangle
\begin{equation}
\label{eq:HKR}
E_*\otimes N^\vee[1]\ra \tau^{\geq -1}\bi_D^*\bi_{D,*}E_*\ra E_*\ra E_*\otimes N^\vee[2].\end{equation}
A priori the HKR class is an element of 
\[\Ext^2_{(X,D|_X,\sA_*|_X)}(E_*,E_*\otimes N^\vee)\]
 which is the obstruction of lifting the splitting module $E_*$ to the first infinitesimal neighborhood of the embedding $\bi_D$.  
\paragraph By Morita equivalence the extension group above is isomorphic to
\[Ext^2_{(X,D|_X)}(\cO_{(X,D|_X)},\cO_{(X,D|_X)}\otimes N^\vee)=H^2((X,D|_X),\cO_{(X,D|_X)}\otimes N^\vee).\]
It is easy to see that $\Hom_{(X,D|_X)}(\cO_{(X,D|_X)}, E'_*)=\Hom_X(\cO_X,E'_0)$ for any parabolic sheaf $E'_*$.  Therefore, the HKR class $\alpha^D_{E_*}$ can be though of as an element of $H^2(X,N^\vee)$.  This element corresponds to the rightmost map of the triangle given by the triangle in \ref{eq:HKR} when $*=0$, i.e, all objects are considered as ordinary sheaves.  Summarizing the above discussion we obtain the following theorem which is a straightforward generalization of Theorem \ref{thm:Azu}.

\begin{Theorem}
\label{thm:lat}
Assume that both maps $i$ and $i_D$ split to first order.  Then, the following statements are equivalent.
\begin{itemize}
\item[(1)] The dg-schemes $(W,D|_W)$ and $(W,D|_W,\sA|_W)$ are isomorphic over $(X\times X,D|_{X\times X})$.
\item[(2)] There exists an isomorphism of dg-autofunctors of $D(X,D|_X)$
\[i_D^*i_{D,*}(-)\iso i_D'^*i'_{D,*}(-)=(-)\otimes \S(N^\vee[1]).\]
\item[(3)] The associated parabolic line bundle is trivial.
\item[(4)] The morphism $i'_D$ splits to first order.
\item[(5)] The HKR class $\alpha^D_{E_*}$ vanishes.
\end{itemize}
\end{Theorem}

\section{Proof of the main theorems}
\label{sec:mainthm}

In this section we prove our main theorems, Theorem \ref{thm:str}, Theorem \ref{thm:alb}, Theorem \ref{thm:equi} and Corollary \ref{cor:end}.

\paragraph Let $X$ be a smooth scheme over a perfect field $k$ of characteristic $p>\dim X$, and $D$ a reduced normal crossing divisor.  We denote by $D'$ the corresponding reduced normal crossing divisor on the Frobenius twist, $X'$.  We consider the Frobenius twist, $X'$ embedded into the vector bundle $T^*X'(\log D')$ as the zero section.  Recall that the parabolic sheaf of crystalline logarithmic differential operators $D(\log D)_*$ can be regarded as a parabolic sheaf of algebras over $T^*X'(\log D')$.  Moreover, $D(\log D)$ is an Azumaya algebra over the logarithmic scheme $(T^*X'(\log D',\pi^*D')$, so that $D(\log D)|_{X'}$ is a split Azumaya algebra.  

\paragraph
We are in the context described in Paragraph \ref{par:setup}.  We compare the derived self-intersection $(W,D|_W)$ corresponding to the embedding 
\[i:(X',D')\ra (T^*X'(\log D'),\pi^*D')\]
and the derived self-intersection corresponding to 
\[i':(X',D')\ra (T^*X'(\log D'),\pi^*D',D(\log D)).\]
We denote the latter space by $(\bW,D'|_W)$, and the map 
\[(X',D',D(\log D)|_{X'})\ra (T^*X'(\log D'),\pi^*D',D(\log D))\]
by $\bi$.  As an easy consequence of Proposition \ref{prop:for} we obtain that the structure sheaf of $(W,D|_W)$ is a formal parabolic sheaf.

\begin{Theorem}
\label{thm:str}
The structure sheaf $\cO_{(W,D|_W)}$ over $X'$ is isomorphic to the dual of the formal complex $\S(\Omega^1_{X'}(\log D')[-1])$ equipped with the trivial parabolic structure.
\end{Theorem}
\medskip

\begin{Proof}
The structure sheaf of $(W,D|_W)$ over $X'$ is given by the object $i^*i_*\cO_{(X',D')}$.  By Remark \ref{rem:ass} all the assumption of Proposition \ref{prop:for} are satisfied for the embedding $i$ implying the statement above.\qed
\end{Proof}
\medskip

Next we compute the associated parabolic line bundle $\sL_*$ of the derived self-intersection corresponding to 
\[i':(X',D')\ra (T^*X'(\log D'),\pi^*D',D(\log D)).\]

\begin{Theorem}
\label{thm:alb}
The associated line bundle $\sL$ is isomorphic to the dual of $F_*\Omega^\sbt_X(\log D)_*$.
\end{Theorem}
\medskip

\begin{Proof}
We have the following sequence of maps
\begin{align*}
F_*\Omega^\sbt_X(\log D)_*&=\RsHom_{(X',D',F_*D_X(\log D)_*)}(F_*\cO_{(X,D)},F_*\cO_{(X,D)})=\\
&=\RsHom_{(X',D',\pi_*D(\log D)_*)}(F_*\cO_{(X,D)},F_*\cO_{(X,D)})=\\
&=\RsHom_{(X',D',\pi_*D(\log D)_*)}(\pi_*\bi_*F_*\cO_{(X,D)},\pi_*\bi_*F_*\cO_{(X,D)})=\\
&=\pi_*\RsHom_{(T^*X'(\log D'),\pi^*D',D(\log D)_*)}(\bi_*F_*\cO_{(X,D)},\bi_*F_*\cO_{(X,D)})\\
\end{align*}
where the first isomorphism is the Koszul duality between $\Omega^\sbt_X(\log D)$ and $D_X(\log D)$, the second is the isomorphism $F_*D_X(\log D)_*=\pi_*D(\log D)_*$ for the bundle map $\pi:T^*X'(\log D')\ra X'$, the third is the identity $\pi\circ i=\id$ and the last one is the conseqence of that $\pi$ is affine.  The map $\bi_*$ has a right adjoint, which we denote by $\bi^!$ (see \cite{Yok} and \cite{AriCalHab2}).  We have
\[F_*\Omega^\sbt_X(\log D)_*=\pi_*i_*\RsHom_{(X',D',D(\log D)_*|_{X'})}(F_*\cO_{(X,D)},\bi^!\bi_*F_*\cO_{(X,D)}).\]
We use again that $\pi\circ i=\id$ to obtain
\[F_*\Omega^\sbt_X(\log D)_*=\RsHom_{(X',D',D(\log D)_*|_{X'})}(F_*\cO_{(X,D)},\bi^!\bi_*F_*\cO_{(X,D)}).\]
We remind the reader that there exists a Morita equivalence between $(X',D')$ and $(X', D', D(\log D))$ given by the functors
\[m_*:Coh(X',D')\ra Coh(D(\log D)|_{(X',D')}):\]
\[M\mapsto M\otimes F_*\cO_{(X,D)}\]
and
\[m^*:Coh(D(\log D)|_{(X',D')})\ra Coh(X',D'):\]
\[ M\mapsto \sHom_{D(\log D)|_{X'}}(F_*\cO_{(X,D)},M)\]
and therefore
\[F_*\Omega^\sbt_X(\log D)_*=m^*\bi^!\bi_*\cO_{(X',D')}.\]
Moreover, the functor $m^*$ is both left and right adjoint to $m_*$ implying that
\[F_*\Omega^\sbt_X(\log D)_*=m^!\bi^!\bi_*\cO_{X'}=i'^!i'_*\cO_{(X',D')}\]
completing the proof.\qed
\end{Proof}
\medskip

We conclude the paper by applying Theorem \ref{thm:lat} to our situation.

\begin{Theorem}
\label{thm:equi}
Let $X$ be a smooth variety over a perfect field of characteristic $p>\dim X$, with a reduced normal crossing divisor $D$.  Then, the following statements are equivalent.
\begin{itemize}
\item[(1)] The logarithmic scheme $(X,D)$ lifts to $W_2(k)$.
\item[(2)] The class of the extension
\[0\ra \cO_{X'}\ra F_*\cO_X\ra F_*Z^1\ra \Omega^1_{X'}(\log D')\ra 0\]
vanishes.  Here $Z^1$ denotes the image of $d$ inside $F_*\Omega^1_X(\log D)$.
\item[(3)] The map $i'$ splits to first order.
\item[(4)] The associated line bundle is trivial.
\item[(5)] The parabolic sheaf of algebras $D(\log D)_*$ splits on the first infinitesimal neighborhood of $(X',D')$ inside $(T^*X'(\log D'),\pi^*D')$.
\item[(6)] The complex $F_*\Omega^\sbt_X(\log D)_*$ is a formal parabolic sheaf equipped with the trivial parabolic structure.
\item[(7)] There exists an isomorphism
\[F_*\Omega^\sbt_X(\log D)_*\iso \S(\Omega^1_{X'}(\log D')[-1])_*\]
in $D(X',D')$ where the sheaf $\S(\Omega^1_{X'}(\log D')[-1])$ is equipped with the trivial parabolic structure.
\end{itemize}
\end{Theorem}
\medskip

\begin{Proof}
The HKR class corresponding to the embedding $i'$ is the dual of the extension class in 2.  Therefore, the equivalence \[(2)\Leftrightarrow (3)\Leftrightarrow(4)\Leftrightarrow(5)\Leftrightarrow(6)\Leftrightarrow(7)\]
follows from Theorem \ref{thm:lat}, Theorem \ref{thm:str} and Theorem \ref{thm:alb}.  The equivalence $(1)\Leftrightarrow(2)$ is proved in \cite{EsnVie}.
\qed
\end{Proof}
\medskip

\paragraph[Remark:] A particular interesting feature of the theorem above is part (7).  The two parabolic complexes $F_*\Omega^\sbt_X(\log D)_*$ and $\S(\Omega^1_{X'}(\log D')[-1])_*$ have different filtration, the former one has weights $\alpha_i=\frac{i}{p+1}$, the latter has the trivial parabolic structure.  On the other hand under quasi-isomorphism the filtrations may change.  We illustrate the phenomenon in the simplest case.  Let $X=\Spec k[x]$ and $D=\{0\}$.  Then the pushforward of the logarithmic de Rham complex is given by
\[\left(x^pk[x]\inj x^{p-1}k[x]\inj...\inj k[x]\right)\xrightarrow{d}\left(x^p\frac{k[x]dx}{x}\inj...\inj \frac{k[x]dx}{x}\right).\]
The kernel of $d$ is the parabolic sheaf
\[x^pk[x^p]\inj x^pk[x^p]\inj ...\inj x^pk[x^p]\inj k[x^p]\]
which is isomorphic $\cO_{(X',D')}$.  Similar calculations can be done for the cokernel.

\paragraph
As an immediate consequence we obtain that the Hodge-to-de Rham spectral sequence for the logarithmic de Rham complex degenerates at ${}^1E$.

\begin{Corollary}
\label{cor:end}
Let $X$ be a smooth variety over a perfect field of characteristic $p>\dim X$, and $D$ a reduced normal crossing divisor on $X$.  Assume that $(X,D)$ lifts to $W_2(k)$.  Then,
\[R^*\Gamma(X,\Omega^\sbt_X(\log D))=\bigoplus_{p+q=*} H^p(X,\Omega^q_{X}(\log D)).\]
\end{Corollary}
\medskip 

\begin{Proof}
By Theorem \ref{thm:equi} there exists an isomorphism
\[F_*\Omega^\sbt_X(\log D)_*\iso \S(\Omega^1_{X'}(\log D')[-1])_*.\]
In particular setting $*=0$ we obtain an isomorphism
\[F_*\Omega^\sbt_X(\log D)\iso \S(\Omega^1_{X'}(\log D')[-1])\]
of objects in $D(X')$.  Thus,
\[R^*\Gamma(X',F_*\Omega^\sbt_X(\log D))=R^*\Gamma(X',\S(\Omega^1_{X'}(\log D')[-1])=\bigoplus_{p+q=*} H^p(X',\Omega^q_{X'}(\log D')).\]
On one hand, since $F$ is affine we have 
\[R^*\Gamma(X',F_*\Omega^\sbt_X(\log D))=R^*\Gamma(X,\Omega^\sbt_X(\log D)).\]
On the other hand, since $X$ and $X'$ are abstractly isomorphic, we have
\[\bigoplus_{p+q=*} H^p(X',\Omega^q_{X'}(\log D'))=\bigoplus_{p+q=*} H^p(X,\Omega^q_{X}(\log D))\]
completing the proof.\qed
\end{Proof}

\end{document}